# Disjoint $X$-paths in bidirected graphs

Jana K. Nickel

February 27, 2025

### Abstract

Let $B$ be a bidirected multigraph with signing $\sigma$, let $X$ be a set of vertices in $B$, and let $k$ be a non-negative integer. For any pair of vertex sets $S, T \subset V(B)$ satisfying $X \cap S = X \cap T$, we denote by $B_{S,T}$ the multigraph with the same vertex set as $B$ and with edge set consisting of those edges $e$ of $B$ each of whose endvertices $v$ satisfies $v \notin S \cup T$ or $v \in S \setminus T$, $\sigma(v, e) = -$ or $v \in T \setminus S$, $\sigma(v, e) = +$. We prove that $B$ admits a set of $k$ pairwise disjoint $X$-paths if and only if for any $S, T \subseteq V(B)$ with $X \cap S = X \cap T$, the inequality $|S \cap T| + \sum \lfloor \frac{1}{2} |V(C) \cap (X \cup S \cup T)| \rfloor \geq k$ holds where the sum is indexed by the components of $B_{S,T}$. This result is a generalization of a result of Gallai from undirected graphs to bidirected ones. Furthermore, we will deduce from this a kind of an Erdős-Pósa property for $X$-paths in bidirected multigraphs.

## 1 Introduction

Given an undirected graph $G$, a theorem of Gallai [3] states that for every vertex set $X \subseteq V(G)$ and every integer $k \geq 0$, there exist $k$ pairwise disjoint $X$-paths in $G$ if and only if any vertex set $S \subseteq V(G)$ satisfies $|S| + \sum \lfloor \frac{1}{2} |V(C) \cap X| \rfloor \geq k$ where the sum is taken over the set $\mathcal{C}(G-S)$ of components of the subgraph $G - S$ of $G$. An $X$-*path* is a non-trivial path which intersects $X$ precisely in its two endvertices. Gallai's result does not apply to bidirected graphs, no matter whether components are replaced by strong or weak components (defined in Section 2). In the case of strong components, the example on the left-hand side of Figure 1.1 illustrates that the condition on the upper bound of $k$ is not in general necessary for the existence of $k$ pairwise disjoint $X$-paths. Considering weak components, the example on the right-hand side of the figure shows that the condition is not in general sufficient.

If we replace edges leaving a vertex $u$ with negative sign and entering a vertex $v$ with positive sign by directed edges from $u$ to $v$, then the corresponding examples in Figure 1.1 reveal that

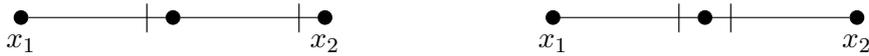

Figure 1.1: On the left-hand side, putting $X \coloneqq \{x_1, x_2\}$, then for the empty set $S$, we have $|S| + \sum \lfloor \frac{1}{2}|V(C) \cap X|\rfloor = 0$ where the sum runs over the set of strong components of the given bidirected graph $B$. However, $B$ contains an $X$-path, so for $k = 1$, the condition on the upper bound of $k$ in Gallai's theorem is not necessary. On the right-hand side, defining $X$ as before, then for any $S \subseteq V(B)$, we have $|S| + \sum \lfloor \frac{1}{2}|V(C) \cap X|\rfloor \geq 1$ where the sum runs over the set of weak components of $B - S$. But $B$ contains no $X$-path, so for $k = 1$, the condition is not sufficient.

Gallai's result does not either hold for directed graphs. In his paper [8], Kriesell found a both necessary and sufficient condition similar to that in Gallai's theorem by refining the notion of vertex separation. Given a digraph $D$ and a vertex set $X \subseteq V(D)$, instead of taking single sets $S$ of vertices in $D$ and the digraph $D - S$, which is obtained from $D$ by deleting the vertices in $S$ together with all edges whose terminal or initial vertex lies in $S$, Kriesell considers pairs $(S, T)$ of vertex sets $S, T \subseteq V(D)$ with $X \cap S = X \cap T$ and the digraph $D - (S, T)$ obtained from $D$ by deleting all edges whose terminal vertex lies in $S$ and whose initial vertex lies in $T$. Concerning the various types of components of a digraph, he chooses the notion of weak components. Concretely, Kriesell proves that for each subset $X \subseteq V(D)$ and each integer $k \geq 0$, there exist $k$ pairwise disjoint $X$-paths if and only if for any subsets $S, T \subseteq V(D)$ satisfying $X \cap S = X \cap T$, we have $|S \cap T| + \sum \lfloor \frac{1}{2}|V(C) \cap (X \cup S \cup T)|\rfloor \geq k$ where the sum is taken over the weak components of the subdigraph $D - (S, T)$ of $D$. The proof is done in two steps. First, Kriesell defines a certain auxiliary digraph $H_D(X)$ obtained from $D$ by splitting the vertices outside $X$ and shows that the existence of $k$ pairwise disjoint $X$-paths in $D$ is equivalent to the existence of a matching in $H_D(X)$ of cardinality $k + |V(D) \setminus X|$. Second, he exploits Tutte's characterization 3.2 of the existence of a matching of a specific cardinality and relates this to the desired characterization of the existence of $k$ pairwise disjoint $X$-paths in $D$.

In the paper at hand, we generalize Kriesell's theorem and proof to bidirected graphs. These objects were first introduced and studied by Kotzig in a series of papers [5, 6, 7] and were rediscovered independently by Edmonds and Johnson [4]. They are defined as (undirected) graphs together with two signs for every edge, one at each endvertex. In particular, every bidirected graph has a unique underlying undirected graph. Moreover, we can regard bidirected graphs as a generalization of directed graphs inasmuch as every directed graph gives rise to a bidirected graph where any directed edge was replaced by an edge with negative sign at the initial vertex and positive sign at the terminal vertex, while bidirected graphs may additionally contain edges with the same sign at both endvertices. We will mostly consider bidirected graphs which are allowed to have multiple edges with the same pair of signs. For the sake of clarity, we refer to them as bidirected multigraphs. However, we do not allow them to contain loops.

For the proof of our generalized Gallai theorem for bidirected multigraphs, Theorem 3.1, we need a substitute for Kriesell's digraph $D - (S, T)$. This digraph can be equivalently described as the spanning subdigraph of $D$ containing those directed edges $uv$ whose initial vertex $u$ satisfies $u \notin S \cup T$ or $u \in S \setminus T$ and whose terminal vertex $v$ satisfies $v \notin S \cup T$ or $v \in T \setminus S$. In the process of adapting Kriesell's proof to bidirected multigraphs, the multigraph $B_{S,T}$ with vertex set $V(B)$ and edge set consisting of those edges $e$ of $B$ each of whose endvertices $v$ satisfies $v \notin S \cup T$ or $v \in S \setminus T$, $\sigma(v,e) = -$ or $v \in T \setminus S$, $\sigma(v,e) = +$ turns out to be suitable. Here, $\sigma$ denotes the *signing* of $B$, assigning to each pair of an edge and one of its endvertices a plus or minus sign. If replacing any directed edge of a digraph $D$ from $u$ to $v$ by an edge with sign $-$ at $u$ and sign $+$ at $v$, then for the resulting bidirected graph $B$, $B_{S,T}$ corresponds precisely to the directed subdigraph $D - (S,T) \subseteq D$.

As a consequence of Theorem 3.1, we will prove that $X$-paths in bidirected multigraphs satisfy a kind of an Erdős-Pósa property, namely that there exists a function $f\colon \mathbb{N}_0 \to \mathbb{N}_0$ such that for any bidirected multigraph $B$, for any vertex set $X \subseteq V(B)$ and any integer $k \geq 0$, $B$ admits $k$ pairwise disjoint $X$-paths or a vertex set $Y \subseteq V(B)$ of cardinality at most $f(k)$ such that $B - Y$ contains no $X$-path. We will verify this for $f(k) := 2k - 2$, see Theorem 4.1. This result was already known for undirectd graphs and first appeared in Gallai's paper [3]. More general, a family $\mathcal{F}$ of graphs is said to satisfy the *Erdős-Pósa property* iff there is a function $f\colon \mathbb{N}_0 \to \mathbb{N}_0$ such that for each integer $k \geq 0$, every graph $G$ contains $k$ pairwise disjoint subgraphs in $\mathcal{F}$ or a vertex set $Y \subseteq V(G)$ of cardinality $|Y| \leq f(k)$ such that $G - Y$ contains no subgraph in $\mathcal{F}$. Any such $f$ is called an *Erdős-Pósa function* for $\mathcal{F}$.

**Acknowledgement.** I would like to thank Paul Knappe for simplifying some arguments in the proof of Theorem 3.1.

## 2 Preliminaries: Bidirected multigraphs

We use standard graph-theoretic terminology from [2]. The following definitions for bidirected multigraphs are based on Section 2 of [9].

A *bidirected multigraph* $B$ is a pair $(G, \sigma)$ consisting of a multigraph $G = (V, E)$ together with a *signing* $\sigma\colon \underline{E}(B) \to \{+, -\}$ where $\underline{E}(B)$ denotes the set of half-edges of $B$, that is, $\underline{E}(B) := \{(v, e) \in V \times E \mid v \text{ is an endvertex of } e\}$. Given a half-edge $(v, e)$ of $B$, we say that $e$ has sign $\sigma(v, e)$ at $v$. Moreover, edges with signs $\alpha$ and $\beta$ at their endvertices are referred to as $(\alpha, \beta)$-edges. The *vertex set* $V(B)$ of $B$ is the vertex set $V = V(G)$ of the *underlying multigraph* $G$ and the *edge set* $E(B)$ of $B$ is the edge set $E = E(G)$ of $G$. We call the elements of $V(B)$ and $E(B)$ the *vertices* and *edges* of $B$, respectively. If $B$ has no multiple edges with the same pair of signs, we refer to $B$ as a *bidirected graph*.

When drawing a bidirected multigraph, we omit any negative sign while representing any positive sign as a bar perpendicular to the respective edge near the respective endvertex. For an example, see Figure 3.1.

Let $B = (G, \sigma)$ be a bidirected multigraph. A *bidirected sub(multi)graph* of $B$ is a bidirected (multi)graph $B' = (G', \sigma')$ where $G'$ is a submultigraph of $G$ and $\sigma'$ is the restriction of $\sigma$ to $\underline{E}(B') \subseteq \underline{E}(B)$. Given a vertex set $Y \subset V(B)$, we denote by $B - Y$ the bidirected submultigraph of $B$ whose underlying multigraph is $G - Y$, that is, the multigraph obtained from $G$ by deleting all vertices in $Y$.

A *path* in $B$ is a path $P = v_0 e_1 v_1 \ldots e_\ell v_\ell$ in the underlying multigraph $G$ such that for each internal vertex $v \in \{v_1, \ldots, v_{\ell-1}\}$, the two half-edges of $P$ containing $v$ have distinct signs, that is, $\sigma(v_i, e_i) \neq \sigma(v_i, e_{i+1})$ for all $i \in \{1, \ldots, \ell-1\}$. We also call $P$ a $(v_0, \sigma(v_0, e_1)) - (v_\ell, \sigma(v_\ell, e_\ell))$ *path*. The *length* of $P$ is the number $\ell$ of its edges. We say that $P$ is *trivial* iff it has length zero. Given a vertex set $X \subseteq V(B)$, an *$X$-path* in $B$ is a non-trivial path $P = v_0 e_1 v_1 \ldots e_\ell v_\ell$ in $B$ whose endvertices $v_0$ and $v_\ell$ lie in $X$ whereas $X$ contains no other vertex of $P$. Note that every path in $B$ defines a bidirected subgraph of $B$ in the obvious way, and we usually do not distinguish between $P$ and the corresponding bidirected graph.

Similarly as for directed multigraphs, there are different notions of connectivity for their bidirected counterparts. A bidirected multigraph $B$ is *strongly connected* iff it is non-empty and for any two vertices $u, v \in V(B)$, there exist signs $\alpha$ and $\beta$ such that $B$ contains both a $(u, \alpha) - (v, \beta)$ path and a $(u, -\alpha) - (v, -\beta)$ path. It is said to be *weakly connected* iff its underlying multigraph is connected (and thus non-empty). Clearly, every strongly connected bidirected multigraph is in particular weakly connected. The *strong (weak) components* of $B$ are its maximal strongly (weakly) connected bidirected submultigraphs. Note that the weak components of $B$ are precisely the components of its underlying multigraph. For an (undirected) multigraph $H$, we denote by $\mathcal{C}(H)$ the set of its components, and for a bidirected multigraph $B = (G, \sigma)$, we write $\mathcal{C}(B) \coloneqq \mathcal{C}(G)$.

## 3 A characterization of the existence of $k$ pairwise disjoint $X$-paths

Our first main theorem, stated below, is a bidirected generalization of Gallai's theorem for undirected graphs [3] and of Kriesell's theorem for directed graphs [8].

**Theorem 3.1.** *Let $B$ be a bidirected multigraph with signing $\sigma \colon \underline{E}(B) \to \{+, -\}$, let $X \subseteq V(B)$ be a vertex set, and let $k$ be a non-negative integer. Then $B$ admits a set of $k$ pairwise disjoint $X$-paths if and only if for any pair of vertex sets $S, T \subseteq V(B)$ with $X \cap S = X \cap T$, we have*
$$|S \cap T| + \sum_{C \in \mathcal{C}(B_{S,T})} \left\lfloor \tfrac{1}{2} |V(C) \cap (X \cup S \cup T)| \right\rfloor \geq k$$
*where $B_{S,T}$ denotes the undirected multigraph $(V(B), \{e \in E(B) \mid \text{each endvertex } v \text{ of } e \text{ satisfies } v \notin S \cup T \text{ or } v \in S \setminus T, \sigma(v, e) = - \text{ or } v \in T \setminus S, \sigma(v, e) = +\})$.*

The rough structure of our proof is based on [8]. In particular, we will translate the problem of finding pairwise disjoint $X$-paths to a matching problem. In this context, the lemma below will play an important role. It follows from results of Tutte [10] and Berge [1] and characterizes the existence of a matching of a given cardinality in a graph by a simple inequality. We will then relate this to the desired inequality in the theorem.

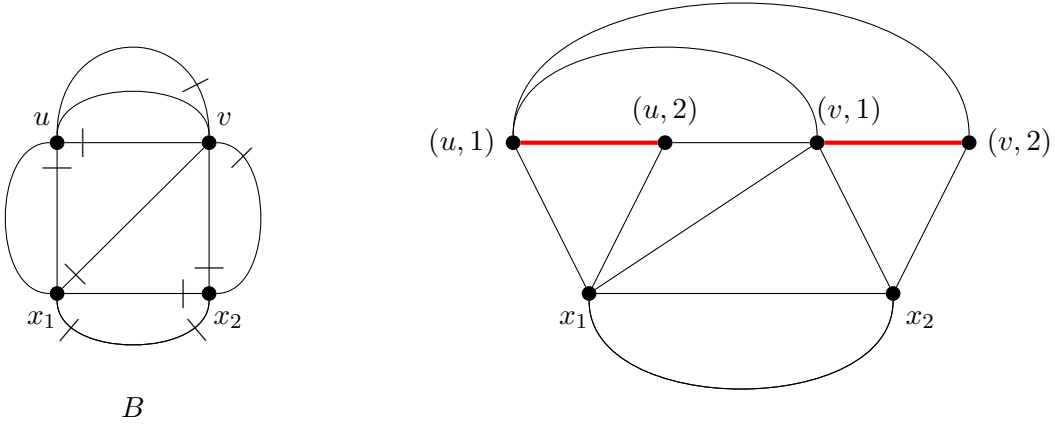

Figure 3.1: Construction of the auxiliary multigraph $H = H_B(X)$ in the proof of Theorem 3.1. In this example, $X$ is the set $\{x_1, x_2\}$, and the matching $M_0$ in $H$ consists of the two red edges.

**Lemma 3.2.** *A graph $H$ has a matching of cardinality $m \in \mathbb{N}_0$ if and only if every vertex set $U \subseteq V(H)$ satisfies*

$$|U| + \sum_{C \in \mathcal{C}(H-U)} \left\lfloor \tfrac{1}{2} |C| \right\rfloor \geq m.$$

This result obviously extends to multigraphs.

*Proof of Theorem 3.1.* We define an auxiliary (undirected) multigraph $H = H_B(X)$ as follows. Its vertex set is $V(H) := X \cup ((V(B) \setminus X) \times \{1, 2\})$. Let $p \colon V(B) \times \{1, 2\} \to V(H)$ be the surjective map given by

$$p(v, i) := \begin{cases} v & v \in X \\ (v, i) & v \notin X. \end{cases}$$

For any vertex $v \in V(B) \setminus X$, let $\widetilde{e}^v$ be an edge in $H$ between $(v, 1)$ and $(v, 2)$, and for any edge $e \in E(B)$ with endvertices $u$ and $v$, let $\widetilde{e}$ be an edge in $H$ with endvertices $p(u, i_{u,e})$ and $p(v, i_{v,e})$ where the index $i_{v,e} \in \{1, 2\}$ is defined as

$$i_{v,e} := \begin{cases} 1 & \sigma(v, e) = - \\ 2 & \sigma(v, e) = +. \end{cases}$$

We put

$$E(H) := \{\widetilde{e}^v \mid v \in V(B) \setminus X\} \cup \{\widetilde{e} \mid e \in E(B)\}.$$

In other words, $H$ is the multigraph obtained from the underlying multigraph of $B$ by replacing each vertex $v$ outside $X$ by two new vertices $(v, 1)$ and $(v, 2)$ together with an edge $\widetilde{e}^v$ between $(v, 1)$ and $(v, 2)$, and by replacing each edge $e$ with endvertices $u$ and $v$ by an edge $\widetilde{e}$ between $p(u, i_{u,e})$ and $p(v, i_{v,e})$. By construction, the set $M_0 := \{\widetilde{e}^v \mid v \in V(B) \setminus X\} \subset E(H)$ is a matching in $H$ An example of $H$ together with the matching $M_0$ is shown in Figure 3.1.

Observe that the $X$-paths in the bidirected multigraph $B$ are in bijection with the $M_0$-alternating $X$-paths in the undirected multigraph $H$, that is, the $X$-paths in $H$ whose edges alternate between $M_0$ and $E(H) \setminus M_0$. In fact, an $X$-path $P = v_0 e_1 v_1 \ldots v_{\ell-1} e_\ell v_\ell$ in $B$ corresponds to the $M_0$-alternating $X$-path

$$\theta(P) := v_0 \widetilde{e}_1(v_1, i_{v_1, e_1}) \widetilde{e}^{v_1}(v_1, i_{v_2, e_1}) \widetilde{e}_2(v_2, i_{v_2, e_2}) \widetilde{e}^{v_2}(v_2, i_{v_3, e_2}) \ldots (v_{\ell-1}, i_{v_\ell, e_{\ell-1}}) \widetilde{e}_\ell v_\ell$$

in $H$. Furthermore, note that two $X$-paths $P$ and $P'$ in $B$ are disjoint if and only if the two $X$-paths $\theta(P)$ and $\theta(P')$ in $H$ are disjoint.

**Claim 3.3.** *For every integer $k \geq 0$, $B$ contains $k$ pairwise disjoint $X$-paths if and only if $H$ admits a matching of cardinality $k + |V(B) \setminus X|$.*

*Proof.* The matching $M_0 = \{\widetilde{e}^v \mid v \in V(B) \setminus X\}$ in $H$ has cardinality $|V(B) \setminus X|$, so for $k = 0$, both conditions hold. Let $k \geq 1$, and let $\mathcal{P}$ be a set of $k$ pairwise disjoint $X$-paths in $B$. For every $P \in \mathcal{P}$, none of the two endvertices of the $M_0$-alternating $X$-path $\theta(P)$ in $H$ is matched by $M_0$. Thus, the symmetric difference of $M_0$ and $P$ is a matching in $H$ of cardinality $|M_0|+1$. Starting with $M_0$, we successively take the symmetric difference with $\theta(P)$ for each $P \in \mathcal{P}$. Since the paths $\theta(P)$ with $P \in \mathcal{P}$ are pairwise disjoint, this provides a matching in $H$ of cardinality $|\mathcal{P}| + |M_0| = k + |V(B) \setminus X|$.

Conversely, let $M$ be a matching in $H$ of cardinality $k + |V(B) \setminus X|$. We denote by $H'$ the submultigraph $(V(H), M_0 \cup M) \subset H$. Then every component of $H'$ is a path or an even cycle in $H$ whose edges alternate between $M_0 \setminus M$ and $M \setminus M_0$. Let us show that at least $k$ of them are $X$-paths in $H$. Clearly, every even cycle in $H'$ has the same number of edges in $M_0$ and in $M$. We write $\mathcal{Q} \subseteq \mathcal{C}(H')$ for the components of $H'$ that are $X$-paths in $H$. By virtue of the bijection between the $M_0$-alternating $X$-paths in $H$ and the $X$-paths in $B$, we see that every $Q \in \mathcal{Q}$ has precisely one edge more in $M$ than in $M_0$ whereas for any path $Q \in \mathcal{C}(H')$ that is not an $X$-path in $H$, the number of edges of $Q$ in $M$ is less than or equal to the number of its edges in $M_0$. Using the fact that

$$M_0 \cup M = E(H') = \bigcup_{Z \in \mathcal{C}(H')} E(Z) = \bigcup_{\text{cycles } C \in \mathcal{C}(H')} E(C) \cup \bigcup_{Q \in \mathcal{Q}} E(Q) \cup \bigcup_{\text{paths } Q \in \mathcal{C}(H') \setminus \mathcal{Q}} E(Q),$$

we therefore have

$$\begin{aligned} |M| &= \sum_{\text{cycles } C \in \mathcal{C}(H')} |E(C) \cap M| + \sum_{Q \in \mathcal{Q}} |E(Q) \cap M| + \sum_{\text{paths } Q \in \mathcal{C}(H') \setminus \mathcal{Q}} |E(Q) \cap M| \\ &\leq \sum_{\text{cycles } C \in \mathcal{C}(H')} |E(C) \cap M_0| + \sum_{Q \in \mathcal{Q}} (1 + |E(Q) \cap M_0|) + \sum_{\text{paths } Q \in \mathcal{C}(H') \setminus \mathcal{Q}} |E(Q) \cap M_0| \\ &= |\mathcal{Q}| + |M_0|, \end{aligned}$$

and hence, $|\mathcal{Q}| \geq |M| - |M_0| = (k + |V(B) \setminus X|) - |V(B) \setminus X| = k$. This shows that $\mathcal{C}(H')$ contains at least $k$ $X$-paths in $H$. By construction, they alternate between $M_0 \setminus M$ and

$M \setminus M_0 \subset E(H) \setminus M_0$, and as components of the submultigraph $H' \subseteq H$, these paths are pairwise disjoint, so we deduce that $H$ admits $k$ pairwise disjoint $M_0$-alternating $X$-paths, which means that the bidirected multigraph $B$ admits $k$ pairwise disjoint $X$-paths. △

By Lemma 3.2, the existence of a matching in $H$ of cardinality $k + |V(B) \setminus X|$ is in turn equivalent to the requirement that for any vertex set $U \subseteq V(H)$, the inequality

$$|U| + \sum_{C \in \mathcal{C}(H-U)} \left\lfloor \tfrac{1}{2} |C| \right\rfloor \geq k + |V(B) \setminus X| \tag{3.4}$$

holds. Given $S, T \subseteq V(B)$ with $X \cap S = X \cap T$, we define a map $\gamma = \gamma_{S,T} \colon V(B) \to 2^{V(H)}$ by

$$\gamma(v) := p(\{v\} \setminus T \times \{1\}) \cup p(\{v\} \setminus S \times \{2\}) = \begin{cases} \varnothing & v \in S \cap T \\ \{(v,1)\} & v \in S \setminus T \\ \{(v,2)\} & v \in T \setminus S \\ \{v\} & v \in X \setminus (S \cup T) \\ \{(v,1),(v,2)\} & v \in V(B) \setminus (X \cup S \cup T). \end{cases}$$

Consider the set $U \coloneqq p(T \times \{1\}) \cup p(S \times \{2\}) \subseteq V(H)$. We have $U = (X \cap S \cap T) \sqcup (T \setminus X \times \{1\}) \sqcup (S \setminus X \times \{2\})$ and therefore

$$|U| = |X \cap S \cap T| + |T \setminus X| + |S \setminus X|. \tag{3.5}$$

Moreover, note that $\bigcup_{v \in V(B)} \gamma(v) = V(H) \setminus U$. Let $B_{S,T}$ be the (undirected) multigraph with vertex set $V(B)$ and edge set consisting of those edges $e$ of $B$ each of whose endvertices $v$ satisfies $v \notin S \cup T$ or $v \in S \setminus T$, $\sigma(v,e) = -$ or $v \in T \setminus S$, $\sigma(v,e) = +$. Then the following conditions hold:

(1) For any edge $e \in E(B_{S,T})$ with endvertices $u$ and $v$ and for any $a \in \gamma(u)$ and $b \in \gamma(v)$, there exists an $a$-$b$ path in $B_{S,T}$.

(2) For any $u, v \in V(B)$ and for any edge $f \in E(H)$ with endvertices $a \in \gamma(u)$ and $b \in \gamma(v)$, there exists a $u$-$v$ path in $B_{S,T}$.

We denote by $\mathcal{C}'(B_{S,T})$ the set of components of the multigraph $B_{S,T}$ that do not consist of a single vertex in $S \cap T$. Using (1), one can show that for every $C \in \mathcal{C}'(B_{S,T})$, there exists some $D \in \mathcal{C}(H - U)$ with $\bigcup_{v \in V(C)} \gamma(v) \subseteq V(D)$, and as a consequence of (2), for every $D \in \mathcal{C}(H - U)$, there exists some $C \in \mathcal{C}'(B_{S,T})$ with $V(D) \subseteq \bigcup_{v \in V(C)} \gamma(v)$. This yields the equality

$$\{ \bigcup_{v \in V(C)} \gamma(v) \mid C \in \mathcal{C}'(B_{S,T}) \} = \{ V(D) \mid D \in \mathcal{C}(H - U) \}. \tag{3.6}$$

Given any $C\in\mathcal{C}'(B_{S,T})$, we observe that

$$\bigcup_{v\in V(C)} \gamma(v) = p(V(C)\setminus T\times\{1\}) \cup p(V(C)\setminus S\times\{2\})$$
$$= (X\cap V(C)\setminus (S\cup T))\sqcup (V(C)\setminus (X\cup T)\times\{1\})\sqcup (V(C)\setminus (X\cup S)\times\{2\}),$$

and therefore,

$$\Big|\bigcup_{v\in V(C)}\gamma(v)\Big| = |(X\cap V(C)\setminus (S\cup T)| + |V(C)\setminus (X\cup T)| + |V(C)\setminus (X\cup S)|$$
$$= |(X\cap V(C)\setminus (S\cup T)| + |S\cap V(C)\setminus (X\cup T)| + |T\cap V(C)\setminus (X\cup S)|$$
$$\quad + 2|V(C)\setminus (X\cup S\cup T)|$$
$$= |V(C)\cap (X\cup S\cup T)| + 2|V(C)\setminus (X\cup S\cup T)|. \qquad (3.7)$$

We obtain

$$|U| - |V(B)\setminus X| + \sum_{D\in\mathcal{C}(H-U)} \lfloor\tfrac{1}{2}|D|\rfloor$$
$$\stackrel{(3.6)}{=} |U| - |V(B)\setminus X| + \sum_{C\in\mathcal{C}'(B_{S,T})} \Big\lfloor\tfrac{1}{2}\Big|\bigcup_{v\in V(C)}\gamma(v)\Big|\Big\rfloor$$
$$\stackrel{(3.7)}{=} |U| - |V(B)\setminus X| + \sum_{C\in\mathcal{C}'(B_{S,T})}\lfloor\tfrac{1}{2}|V(C)\cap (X\cup S\cup T)|\rfloor + \sum_{C\in\mathcal{C}'(B_{S,T})}|V(C)\setminus (X\cup S\cup T)|$$
$$\stackrel{(3.5)}{=} |X\cap S\cap T| + |S\setminus X| + |T\setminus X| - |V(B)\setminus X| + \sum_{C\in\mathcal{C}(B_{S,T})} |V(C)\setminus (X\cup S\cup T)|$$
$$\quad + \sum_{C\in\mathcal{C}(B_{S,T})} \lfloor\tfrac{1}{2}|V(C)\cap (X\cup S\cup T)|\rfloor$$
$$= |X\cap S\cap T| + |S\setminus X| + |T\setminus X| + |V(B)\setminus (X\cup S\cup T)| - |V(B)\setminus X|$$
$$\quad + \sum_{C\in\mathcal{C}(B_{S,T})} \lfloor\tfrac{1}{2}|V(C)\cap (X\cup S\cup T)|\rfloor$$
$$= |S\cap T| + \sum_{C\in\mathcal{C}(B_{S,T})} \lfloor\tfrac{1}{2}|V(C)\cap (X\cup S\cup T)|\rfloor. \qquad (3.8)$$

By virtue of this, if every vertex set $U\subseteq V(H)$ satisfies the inequality in (3.4), then for any vertex sets $S,T\subseteq V(B)$ with $X\cap S = X\cap T$, taking $U:= p(T\times\{1\})\cup p(S\times\{2\})$ as above, then $|S\cap T| + \lfloor\tfrac{1}{2}|V(C)\cap (X\cup S\cup T)|\rfloor \geq k$. Conversely, if all $S,T\subseteq V(B)$ with $X\cap S = X\cap T$ satisfy this latter inequality, then for any $U\subseteq V(H)$, defining $S:= \{v\in V(B)\mid p(v,1)\in U\}$ and $T:= \{v\in V(B)\mid p(v,2)\in U\}$, then $S\cap X = U\cap X = X\cap T$ and $U = p(T\times\{1\})\times p(S\times\{2\})$, so by (3.8), $U$ satisfies (3.4).

Together with Claim 3.3 and our observation in (3.4), following from Lemma 3.2, this implies the desired equivalence in Theorem 3.1. $\square$

## 4 The Erdős-Pósa property for $X$-paths in bidirected multigraphs

Here is our second main theorem, saying that $X$-paths in bidirected multigraphs satisfy a kind of an Erdős-Pósa property with Erdős-Pósa function given by $f(k) := 2k - 2$.

**Theorem 4.1.** *Let $B$ be a bidirected multigraph, let $X \subseteq V(B)$ be a vertex set, and let $k \geq 0$ be an integer. Then there exist $k$ pairwise disjoint $X$-paths in $B$ or a vertex set $Y \subseteq V(B)$ of cardinality $|Y| \leq 2k - 2$ such that the bidirected multigraph $B - Y$ contains no $X$-path.*

*Proof.* We suppose that $B$ does not admit $k$ pairwise disjoint $X$-paths. Then by Theorem 3.1, there exist vertex sets $S, T \subseteq V(B)$ with $X \cap S = X \cap T$ such that

$$|S \cap T| + \sum_{C \in \mathcal{C}(B_{S,T})} \left\lfloor \tfrac{1}{2} |V(C) \cap (X \cup S \cup T)| \right\rfloor \leq k - 1$$

where $B_{S,T}$ denotes the multigraph $(V(B), \{e \in E(B) \mid \text{each endvertex } v \text{ of } e \text{ satisfies } v \notin S \cup T$ or $v \in S \setminus T, \sigma(v, e) = -$ or $v \in T \setminus S, \sigma(v, e) = +\})$. Given any such component $C \in \mathcal{C}(B_{S,T})$, we choose a set $Z(C)$ of all but one elements in $V(C) \cap (X \cup S \cup T)$, unless this set is empty, in which case, we put $Z(C) := \varnothing$. We now define

$$Y := (S \cap T) \cup \bigcup_{C \in \mathcal{C}(B_{S,T})} Z(C) \subseteq V(B).$$

For any $C \in \mathcal{C}(B_{S,T})$, we have $|Z(C)| \leq 2 \lfloor \tfrac{1}{2} |V(C) \cap (X \cup S \cup T)| \rfloor$, which implies that

$$|Y| \leq |S \cap T| + 2 \sum_{C \in \mathcal{C}(B_{S,T})} \left\lfloor \tfrac{1}{2} |V(C) \cap (X \cup S \cup T)| \right\rfloor \leq 2(k-1) - |S \cap T| \leq 2k - 2.$$

We prove that $B - Y$ contains no $X$-path. Since $S \cap T \subseteq Y$, it suffices to show that any $X$-path $P$ in $B - (S \cap T)$ has a vertex in $Y$. Write $P = v_0 e_1 v_1 \ldots v_{\ell-1} e_\ell v_\ell$, so $\ell \geq 1$, $v_0, v_\ell \in X \setminus (S \cap T) = X \setminus (S \cup T)$ and $v_1, \ldots, v_{\ell-1} \in V(B) \setminus (X \cup (S \cap T))$. Let us say that a half-edge $(v, e)$ of $B$ is *appropriate* iff $v \in S \setminus T, \sigma(v, e) = -$ or $v \in T \setminus S, \sigma(v, e) = +$. We distinguish the following cases.

(1.1) If $v_1 \in S \setminus T$ and $\sigma(v_1, e_1) = -$, then $e_1 \in E(B')$, so there exists some component $C \in \mathcal{C}(B_{S,T})$ containing $e_1$. We have $v_0, v_1 \in V(C) \cap (X \cup S \cup T)$, so $v_0$ or $v_1$ lies in $Z(C) \subseteq Y$.

(1.2) If $v_1 \in S \setminus T$ and $\sigma(v_1, e_1) = +$, then $\ell \geq 2$ and $\sigma(v_1, e_2) = -$, so the half-edge $(v_1, e_2)$ is appropriate.

(2.1) Similarly as in (1.1), if $v_1 \in T \setminus S$ and $\sigma(v_1, e_1) = +$, then $e_1$ is contained in some $C \in \mathcal{C}(B_{S,T})$, so $v_0$ or $v_1$ lies in $Z(C) \subseteq Y$.

(2.2) Similarly as in (1.2), if $v_1 \in T \setminus S$ and $\sigma(v_1, e_1) = -$, then $\ell \geq 2$ and $\sigma(v_1, e_2) = +$, so $(v_1, e_2)$ is appropriate.

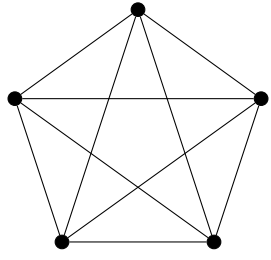

Figure 4.1: A bidirected graph $B$ containing neither $k = 3$ disjoint $V(B)$-paths nor a vertex set $Y \subset V(B)$ with $|Y| < 2k - 2 = 4$ for which $B - y$ has no $V(B)$-path.

(3) If none of the four previous cases occurs, then $v_1 \notin S \cup T$.

For each $i \in \{2, \ldots, \ell\}$ for which $(v_{i-1}, e_i)$ is appropriate or $v_0, \ldots, v_{i-1} \notin S \cup T$, we repeat the above case distinction and inductively continue like this. Assuming that this procedure will not yield any vertex in $V(P) \cap Y$, then all vertices of $P$ necessarily lie in $S \cup T$, so $P$ is entirely contained in $B'$, so $P \subseteq C$ for some $C \in \mathcal{C}(B_{S,T})$. But now, $v_0, v_\ell \in V(C) \cap X$, so $v_0$ or $v_\ell$ lies in $Z(C) \subseteq Y$, a contradiction. We conclude that the path $P$ contains some vertex in $Y$. □

**Remark 4.2.** The upper bound of $2k - 2$ in Theorem 4.1 is tight for every positive integer $k$, that is, for every $k \geq 1$, there exists a bidirected graph $B$ together with some vertex set $X \subseteq V(B)$ such that $B$ neither contains $k$ pairwise disjoint $X$-paths nor a vertex set $Y \subseteq V(B)$ with $|Y| < 2k - 2$ for which $B - Y$ has no $X$-path. For example, consider the bidirected graph $B$ consisting of the complete graph $K^{2k-1}$ on $2k - 1$ vertices with only $(-, -)$-edges (see Figure 4.1 for $k = 3$) and put $X \coloneqq V(B)$. Since every $X$-path has at least two vertices but $|B| = 2k - 1 < 2k$, there are no $k$ pairwise disjoint $X$-paths in $B$. However, for any $Y \subseteq V(B)$ with $|Y| < 2k - 2$, we have $|B - Y| = |B| - |Y| \geq 2$, so $B - Y$ contains an edge, and this is an $X$-path in $B - Y$.

Theorem 3.1 is indeed a generalization of Kriesell's corresponding theorem for directed graphs:

**Corollary 5.1.** *(Kriesell [8]). Let $D$ be a directed multigraph, let $X \subseteq V(D)$, and let $k \geq 0$ be an integer. Then $D$ admits $k$ pairwise disjoint $X$-paths if and only if any $S, T \subseteq V(B)$ with $X \cap S = X \cap T$ satisfy*

$$|S \cap T| + \sum_{C \in \mathcal{C}(D - (E_D^-(S) \cup E_D^+(T)))} \left\lfloor \tfrac{1}{2} |V(C) \cap (X \cup S \cup T)| \right\rfloor \geq k$$

*where $E_D^-(S) \coloneqq E_D(V(D), S)$ is the set of all directed edges of $D$ arriving in $S$ and $E_D^+(T) \coloneqq E_D(T, V(D))$ is the set of all directed edges of $D$ leaving $T$.*

*Proof.* This follows from Theorem 3.1 by considering the bidirected multigraph $B$ that consists of the underlying undirected multigraph of $D$ and the signing $\sigma$ that makes every directed edge of $D$ into a $(-,+)$-edge, that is, for any edge $e \in E(D)$ from $u$ to $v$, we set $\sigma(u,e) \coloneqq -$ and $\sigma(v,e) \coloneqq +$. $\square$

Kriesell shows in [8] that the above corollary in turn implies Gallai's classical result for undirected graphs [3]. Hence, our version for bidirected multigraphs is indeed a generalization of of Gallai's original theorem and of Kriesell's extension to directed graphs. Furthermore, from the Erdős-Pósa property for $X$-paths in bidirected multigraphs, we can easily recover the Erdős-Pósa property for $X$-paths in directed and undirected multigraphs:

**Corollary 5.2.** *Given an undirected or a directed multigraph $G$, a vertex set $X \subseteq V(G)$ and $k \in \mathbb{N}_0$, there exists a set of $k$ pairwise disjoint $X$-paths in $G$ or a vertex set $Y \subseteq V(G)$ with $|Y| \leq 2k - 2$ such that $G - Y$ contains no $X$-path.*

*Proof.* If $G$ is directed, consider $B$ as defined in the proof of the previous corollary and then apply Theorem 4.1. If $G$ is undirected, define $B$ to be the bidirected multigraph obtained from $G$ by replacing each edge with endvertices $u$ and $v$ by two edges $e$ and $f$ where $e$ has sign $-$ at $u$ and sign $+$ at $v$ and $f$ has sign $+$ at $u$ and sign $-$ at $v$. $\square$